\documentclass[12pt,a4paper]{article}
\usepackage[latin1]{inputenc}
\usepackage[english]{babel}
\usepackage{amsmath}
\usepackage{amsfonts}
\usepackage{amssymb}
\usepackage{indentfirst}
\usepackage{dsfont}
\newtheorem{theorem}{Theorem}[section]

\author{Paul Rochet}
\title{Bayesian interpretation of generalized empirical likelihood by maximum entropy}
\date{}

\begin{document}
\maketitle

\begin{abstract} We study a parametric estimation problem related to moment condition models. As an alternative to the generalized empirical likelihood (GEL) and the generalized method of moments (GMM), a Bayesian approach to the problem can be adopted, extending the MEM procedure to parametric moment conditions. We show in particular that a large number of GEL estimators can be interpreted as a maximum entropy solution. Moreover, we provide a more general field of applications by proving the method to be robust to approximate moment conditions.
\end{abstract}

\section{Introduction} 

We consider a parametric estimation problem in a moment condition model. Assume we observe an i.i.d. sample $X_1,...,X_n$ drawn from an unknown probability measure $\mu_0$, we are interested in recovering a parameter $\theta_0 \in \Theta \subset \mathbb R^d$, defined by a set of moment conditions
\begin{equation}\label{eq0gel}   \int \Phi(\theta_0,x) d\mu_0(x) = 0,  \end{equation}
where $\Phi: \theta \times \mathcal X \to \mathbb R^k$ is a known map. This model is involved in many problems in Econometry, notably when dealing with instrumental variables. We refer to \cite{MR888070}, \cite{MR666123}, \cite{MR1272085}, \cite{MR1135146} and \cite{MR2562761}. Two main approaches to the problem have been studied in the literature, namely the generalized method of moments (GMM) and the generalized empirical likelihood (GEL). While the main advantage of GMM relies in its computational feasibility, likelihood-related methods have appeared to be the most efficient in term of small-sample properties. In its original form, the empirical likelihood (EL) of Owen \cite{MR1135146} defines an estimator by a maximum likelihood procedure on a discretized version of the model. As an alternative, GEL replaces the Kullback criterion relative to EL by a $f$-divergence, thus providing a large choice of solutions. A number of estimators corresponding to particular choices of $f$-divergences have emerged in the literature over the last decades, such as the exponential tilting (ET) of Kitamura and Stutzer \cite{MR1458431} and the continuous updating estimator (CUE) of Hansen, Yeaton and Yaron \cite{RePEc:bes:jnlbes:v:14:y:1996:i:3:p:262-80}.

While an attractive feature of GEL is its wide range of solutions, a number of $f$-divergence used in the computation of the GEL estimators are mainly justified by empirical studies and lack a probabilistic interpretation. This issue can be solved by incorporating some prior information to the problem using a Bayesian point of view, as made in \cite{florensrolin}.
In this paper, we investigate a different Bayesian approach to the inverse problem, known as \textit{maximum entropy on the mean} (MEM). Although the method was originally introduced in the frame of exact moment condition models (as opposed here to parametric moment conditions), it appears to provide a natural solution to the problem, expressed as the minimizer of a convex functional on a set of discrete measures and subject to linear constraints. When applied in a particular setting, we show that the MEM approach leads to a GEL solution for which the $f$-divergence is determined by the choice of the prior. As a result, the method gives an alternate point of view on some widely spread estimators such as EL, ET or CUE, as well as a general Bayesian background to GEL. 

In many actual situations, the true moment condition is not exactly known to the statistician and only an approximation is available. It occurs for instance when $\Phi$ has a complicated form that must be evaluated numerically. Simulation-based methods have been implemented to deal with approximate constraints in \cite{MR1803711} and \cite{MR1014539}, in the frame of the generalized method of moments. To our knowledge, the efficiency of GEL in a similar situation has not been studied. In \cite{Devore2}, the MEM procedure is shown to be robust to approximate moment conditions, introducing the approximate maximum entropy on the mean estimator. Seeing GEL as a particular case of MEM, we extend the model in a situation where only an approximation $\Phi_m$ of the true constraint function $\Phi$ is available. We provide sufficient conditions under which the GEL method remains efficient asymptotically when replacing $\Phi$ by its approximation.\bigskip

This paper falls into the following parts. Section \ref{sec2gel} is devoted to the position of the problem. We introduce the maximum entropy method for parametric moment condition models and discuss its close relationship with generalized empirical likelihood in Section \ref{sec3gel}. In Section \ref{GMMsec4}, we discuss the asymptotic efficiency of the method when dealing with an approximate constraint. Proofs are postponed to the Appendix.

\section{Estimation of the parameter}\label{sec2gel}

Let $\mathcal X$ be an open subset of $\mathbb R^q$, endowed with its Borel field $\mathcal B(\mathcal X)$ and let $\mathcal P(\mathcal X)$ denote the set of probability measures on $\mathcal X$. We observe an i.i.d. sample $X_1,...,X_n$ drawn from the unknown distribution $\mu_0$. We want to estimate the parameter $\theta_0 \in \Theta \subset \mathbb R^d$ defined by the moment condition
\begin{equation}\label{eq1} \int_\mathcal X \Phi(\theta_0,x) d\mu_0(x) = 0,  \end{equation}
where $ \Phi: \Theta \times \mathcal X \to \mathbb R^k$ ($k \geq d$) is a known map. To avoid a problem of identifiability, we assume that $\theta_0$ is the unique solution to \eqref{eq1}. This problem has many applications in Econometry, see for instance \cite{MR888070}, \cite{MR666123} and \cite{MR1272085}. The information given by the moment condition \eqref{eq1} can be interpreted to determine the set $\mathcal M$ of possible values for $\mu_0$ (the model). The true value of the parameter being unknown, the distribution of the observations can be any probability measure $\mu$ for which the map $\theta \mapsto \int \Phi(\theta,.) d \mu$ is null for a unique $\theta = \theta(\mu) \in \Theta$. The model is therefore defined as
$$ \mathcal M = \left\{ \mu \in \mathcal P(\mathcal X): \exists ! \ \theta=\theta(\mu) \in \Theta, \textstyle \int \Phi(\theta,.)d \mu =0 \right\},    $$ 
where the map $\mu \mapsto \theta(\mu)$, defined on $\mathcal M$, is the parameter of interest. Let us introduce some notations and assumptions. For $\mu$ a measure and $g$ a function, we shall note $\mu [g] = \int g d \mu$. Let $E$ be an Euclidean space and let $\Vert . \Vert$ denote an Euclidean norm in $E$. For a function $f: \Theta \to E$ and a set $\mathcal S \subseteq \Theta$, we note
$$ \Vert f \Vert_\mathcal S = \sup_{\theta \in \mathcal S} \Vert f(\theta) \Vert.  $$
We assume that the following conditions are fulfilled.
\begin{description}
\item A.1. $\Theta$ is a compact subset of $\mathbb R^d$.
\item A.2. The true value $\theta_0$ of the parameter lies in the interior of $\Theta$.
\item A.3. For all $x \in \mathcal X$, $\theta \mapsto \Phi(\theta,x)$ is continuous on $\Theta$ and the map $x \mapsto \Vert \Phi(.,x) \Vert_\Theta$ is dominated by a $\mu_0$-integrable function.
\item A.4. For all $x \in \mathcal X$, $\theta \mapsto \Phi(\theta,x)$ is twice continuously differentiable in a neighborhood $\mathcal N$ of $\theta_0$ and we note $\nabla \Phi(\theta,x) = \partial \Phi(\theta,x)/\partial \theta \in \mathbb R^{d \times k}$ and $\Psi(\theta,x) = \partial^2 \Phi(\theta,x)/\partial \theta \partial \theta^t \in \mathbb R^{d \times d \times k}$ (where $a^t $ stands for the transpose of $a$). Moreover, we assume that $x \mapsto \Vert \nabla \Phi(.,x) \Vert_\mathcal N $ and $x \mapsto \Vert \Psi(.,x) \Vert_\mathcal N $ are dominated by a $\mu_0$-integrable function.
\item A.5. The matrices 
$$ D := \int_\mathcal X \nabla \Phi(\theta_0,.) d \mu_0 \in \mathbb R^{d\times k} \ \text{ and } \ V:=\int_\mathcal X \Phi(\theta_0,.) \Phi^t(\theta_0,) d \mu_0 \in \mathbb R^{k \times k}$$ 
are of full rank.
\end{description} 
Some issues for estimating $\theta_0$ may be due to the indirect definition of the parameter and these assumptions ensure that the map $\theta(.)$ is sufficiently smooth in a neighborhood of $\mu_0$ for the total variation topology, which will make the asymptotic properties of the GEL estimator easily tractable (see for instance \cite{MR2031017}).

\subsection{Generalized empirical likelihood}
Generalized empirical likelihood (GEL) was first applied to this problem in \cite{MR1272085}, generalizing an idea of \cite{MR1135146}. An estimate $\hat \mu$ of $\mu$ is obtained as an entropic projection of the empirical measure $\mathbb P_n$ onto the model $\mathcal M$. Precisely, for two probability measures $\mu$ and $\nu$ and $f$ a convex function such that $f(1)=f'(1)=0$, define
$$ \mathcal D_f(\nu \vert \mu ) = \int_\mathcal{X} f  \left( \frac{d\nu }{d\mu} \right)  d\mu \ \ \text{if } \nu \ll \mu, \ \ \ \mathcal D_f(\nu \vert \mu ) =  +\infty \: \ \text{otherwise}. $$
Moreover, we define for $\mathcal A \subset \mathcal P (\mathcal X)$, $\mathcal D_f( \mathcal A \vert \mu ) = \inf_{\nu \in \mathcal A} \mathcal D_f(\nu \vert \mu )$. The GEL estimator $\hat \mu$ of $\mu_0$ is the element of the model that minimizes a given $f$-divergence $\mathcal D_f(.,\mathbb P_n)$ with respect to the empirical distribution. Noticing that $\mathcal M = \cup_{\theta \in \Theta} \mathcal M_\theta$ where $ \mathcal M_\theta := \left\{ \mu \in \mathcal P(\mathcal X): \mu[\Phi(\theta,.)] = 0 \right\}$, the GEL estimator $\hat \theta = \theta(\hat \mu)$ follows by
$$ \hat \theta = \arg \min_{\theta \in \Theta} \ \mathcal D_f(\mathcal M_\theta,\mathbb P_n ). $$
Since the set of discrete measures in $\mathcal M_\theta$ is closed and convex, the entropy $\mathcal D_f(\mathcal M_\theta, \mathbb P_n)$ is reached for a unique measure $\hat \mu(\theta)$ in $\mathcal M_\theta$, provided that $\mathcal D_f(\mathcal M_\theta, \mathbb P_n)$ is finite. Then, it appears that computing the GEL estimator involves a two-step procedure. First, build for each $\theta \in \Theta$, the entropic projection $\hat \mu(\theta)$ of $\mathbb P_n$ onto $\mathcal M_\theta$. Then, minimize $\mathcal D_f(\hat \mu(\theta), \mathbb P_n)$ with respect to $\theta$. Since $\hat \mu(\theta)$ is absolutely continuous w.r.t. $\mathbb P_n$ by construction, minimizing $\mathcal D_f(.,\mathbb P_n)$ reduces to finding the proper weights $p_1,...,p_n$ to allocate to the observations $X_1,...,X_n$. This turns into a finite dimensional problem, which can be solved by classical convex optimization tools (see for instance \cite{Kitamura06empiricallikelihood}). In fact, the GEL estimator $\hat \theta$ can be expressed as the solution to the saddle point problem
$$ \hat \theta = \arg \min_{\theta \in \Theta} \ \underset{(\gamma, \lambda) \in \mathbb R \times \mathbb R^{k}}{\text{sup}} \ \gamma - \mathbb P_n \left[ f^*(\gamma + \lambda^t \Phi(\theta,.))\right],$$
where $f^*(x) =  \sup_{y} \left\{ xy - f(y) \right\}$ denotes the convex conjugate of $f$. \bigskip

Note that if the choice of the $f$-divergence plays a key role in the construction of the estimator, it has no influence on its asymptotic efficiency. Indeed, it is shown in \cite{MR1272085} that all GEL estimators are asymptotically efficient, regardless of the $f$-divergence used for their computation. Nevertheless, some situations justify the use of specific $f$-divergences. The empirical likelihood estimator introduced by Owen in \cite{MR1135146} uses the Kullback entropy $\mathcal K(.,.)$ as $f$-divergence, pointing out that minimizing $\mathcal K(.,\mathbb P_n)$ reduces to maximizing likelihood among multinomial distributions. Newey and Smith \cite{MR2031017} remark that a quadratic $f$-divergence leads to the CUE estimator of Hansen Heaton and Yaron \cite{RePEc:bes:jnlbes:v:14:y:1996:i:3:p:262-80}. 

\subsection{Maximum entropy on the mean}\label{sec3gel}

In this section, we study a Bayesian approach to the inverse problem, known as maximum entropy on the mean (MEM) \cite{MR1429928}. The method was developed to estimate a measure $\mu_0$ based the observation of some of its moments. In this framework, it turns out that the MEM estimator of $\mu_0$ can be used to estimate efficiently the parameter $\theta_0$. We shall briefly recall the MEM procedure. Consider an estimator of $\mu_0$ in the form of a weighted version of the empirical measure $\mathbb P_n$,
$$ \mathbb P_n(w) = \frac 1 n \sum_{i=1}^n w_i \ \delta_{X_i},  $$
for $w=(w_1,...,w_n)' \in \mathbb R^n$ a collection of weights. Then, fix a prior distribution $\nu_0$ on the vector of weight $w$ so that each solution $\mathbb P_n(w)$ can be viewed as a realization of the random measure $\mathbb P_n(W)$, where $W$ is drawn from $\nu_0$. This setting enables to incorporate some prior knowledge on the shape or support of $\mu_0$ through the choice of the prior $\nu_0$, as discussed in \cite{MR1429928}. Here, the observations $X_1,...,X_n$ are considered fixed. Actually, it is the moment condition that is used to built the estimator \textit{a posteriori}. In this framework where the true value $\theta_0$ of the parameter is unknown, the information provided by the moment condition reduces to the statement $\mu_0 \in \mathcal M$. So, in order to take this information into consideration, the underlying idea of MEM is to build the estimator $\hat \mu$ as the expectation of $\mathbb P_n(W)$ conditionally to the event $\{\mathbb P_n(W) \in \mathcal M \}$. However, we may encounter some difficulties if this conditional expectation is not properly defined. To deal with this issue, the MEM method replaces the possibly ill-defined conditional expectation by a well-defined estimator, whose construction is motivated by large deviation principles. Precisely, construct the \textit{posterior} distribution $\nu^*$ as the entropic projection of $\nu_0$ onto the set
$$ \Pi(\mathcal M) = \left\{ \mu \in \mathcal P(\mathbb R^n), \ \mathbb E_\mu \left[ \mathbb P_n(W)\right] \in \mathcal M \right\},   $$
where $\mathbb E_\mu \left[ \mathbb P_n(W)\right]$ denotes the expectation of $\mathbb P_n(W)$ when $W$ has distribution $\mu$. The MEM solution to the inverse problem is defined as the expectation of $\mathbb P_n(W)$ under the posterior distribution $\nu^*$,
$$ \hat \mu = \mathbb E_{\nu^*} \left[ \mathbb P_n(W)\right] = \mathbb P_n(\mathbb E_{\nu^*}(W)).   $$
This construction is justified by the large deviation principle stated in Theorem 2.3 in \cite{MR1429928}, which establishes the asymptotic equivalence between $\hat \mu$ and the conditional expectation $\mathbb E_{\nu_0} (\mathbb P_n(W) \vert \ \mathbb P_n(W) \in \mathcal M ) $, whenever it is well defined. The existence of the MEM estimator requires the problem to be \textit{feasible} in the sense that there exists at least one solution $\delta$ in the interior of the convex hull of the support of $\nu_0$, such that $\mathbb P_n(\delta) \in \mathcal M$. This assumption warrants that the set $\Pi(\mathcal M)$ is non-empty and therefore allows the construction of the posterior distribution $\nu^*$. \bigskip

The MEM estimator $ \hat \mu$ lies in the model $\mathcal M$ by construction. As a result, there exists a solution $\hat \theta$ to the moment condition $\hat \mu [\Phi(\theta,.)]=0$. this solution is precisely the MEM estimator of $\theta_0$. In Theorem \ref{gelmemth} below, we give an explicit expression for the MEM estimator $\hat \theta$. We note $\mathds 1 = (1,...,1)^t \in \mathbb R^n$, $\Phi(\theta,X) = (\Phi(\theta,X_1),...,\Phi(\theta,X_n))^t \in \mathbb R^{n \times k}$ and as previously, $\Lambda_\nu$ denotes the log-Laplace transform of $\nu$.
 
\begin{theorem}\label{gelmemth} If the problem is feasible, the MEM estimator $\hat \theta$ is given by
$$ \hat \theta = \arg \min_{\theta \in \Theta} \ \sup_{(\gamma,\lambda) \in \mathbb R \times \mathbb R^{k}}  \left\{ n \gamma - \Lambda_{\nu_0}(\gamma \mathds 1 + \Phi(\theta,X) \lambda )\right\}.$$
In particular, if $\nu_0$ has equal orthogonal marginals, i.e. $\nu_0=\nu^{\otimes n}$ for some probability measure $\nu$ on $\mathbb R$, then
$$ \hat \theta = \arg \min_{\theta \in \Theta} \ \sup_{(\gamma,\lambda) \in \mathbb R \times \mathbb R^{k}} \left\{ \gamma - \mathbb P_n \left[ \Lambda_{\nu}(\gamma + \lambda^t \Phi(\theta,.))\right]\right\}.$$ 
\end{theorem}
The MEM estimator $\hat \theta$ can be expressed as the solution to a saddle point problem, specific to generalized empirical likelihood. Actually, this result points out that maximum entropy on the mean with a particular form of prior $\nu_0 = \nu^{\otimes n}$ leads to a GEL procedure, for which the criterion is the log-Laplace transform of $\nu$. This approach provides a general Bayesian interpretation of GEL. Regularity conditions on the criterion $\Lambda_\nu$ in the GEL framework are reflected through conditions on the prior $\nu$. Indeed, the usual normalization conditions $\Lambda'_\nu(0) = \Lambda''_\nu(0) = 1$ corresponds to taking a prior $\nu$ with mean and variance equal to one, while the normalization $\Lambda_\nu(0)=0$ is imposed by the condition $\nu \in \mathcal P(\mathbb R)$. \bigskip

An interesting choice of the prior is the exponential distribution $d\nu(x) = e^{-x} dx$ for $x>0$. Indeed, observe that if the $W_i$ are i.i.d. with exponential distribution, the likelihood of $\mathbb P_n(W)$ is constant over the set of probability discrete measures $\{ \mathbb P_n(w): \sum^n_{i=1} w_i = n \}$. Hence, an exponential prior can be roughly interpreted as a non-informative prior in this framework. The discrepancy associated to this prior is $\Lambda_\nu(s)= - \log(1-s), \ s < 1$, which corresponds to the empirical likelihood estimator of Owen \cite{MR1135146}. 

The MEM approach also provides a new probabilistic interpretation of some commonly used specific GEL estimators. The exponential tilting of Kitamura and Stutzer \cite{MR1458431} is obtained for a Poisson prior of parameter $1$, for which we have $\Lambda_\nu(s) = e^s -1$. Another example is the Gaussian prior $\nu \sim \mathcal N(1,1)$, leading to the continuous updating estimator of Hansen, Yeaton and Yaron \cite{RePEc:bes:jnlbes:v:14:y:1996:i:3:p:262-80}, as we have in this case $\Lambda_\nu(s) = \frac 1 2 (s-1)^2$. The Gaussian prior allows the discrete measure $\mathbb P_n(W)$ to have negative weights $w_i$ and must be handled with care. Remark however that this is generally not an issue in practice since the solution $\hat \mu$ is implicitly chosen close to the empirical distribution $\mathbb P_n$ and will have all its weights $w_i$ positive with high probability. More examples of classical priors leading to usual discrepancies can be found in \cite{MR1429928}.

\section{Dealing with an approximate operator}\label{GMMsec4}
In many actual applications, only an approximation of the constraint function $\Phi$ is available to the practitioner. This occurs for instance if the moment condition takes a complicated form that can only be evaluated numerically. In \cite{MR1014539}, McFadden suggested a method dealing with approximate constraint in a similar situation, introducing the method of simulated moments (see also \cite{MR1803711}). In \cite{Devore2} and \cite{l2}, the authors study a MEM procedure for linear inverse problems with approximate constraints. Here, we propose to extend the results of \cite{Devore2} and \cite{l2} to the GEL framework, using the connections between GEL and MEM. 
 
We assume that we observe a sequence $\{\Phi_m\}_{m \in \mathbb N}$ of approximate constraints, independent with the original sample $X_1,...,X_n$ and converging toward the true function $\Phi$ at a rate $\varphi_m$. We are interested in exhibiting sufficient conditions on the sequence $\{\Phi_m\}_{m \in \mathbb N}$ under which estimating $\theta_0$ by the GEL procedure remains efficient when the constraint is replaced by its approximation. We discuss the asymptotic properties of the resulting estimates in a framework where both indices $n$ and $m$ simultaneously grow to infinity. \bigskip

\noindent The approximate estimator is obtained by the GEL methodology, replacing the constraint function $\Phi$ by its approximation $\Phi_m$,
\begin{equation} \hat \theta_m = \arg \min_{\theta \in \Theta} \ \sup_{(\gamma,\lambda) \in \mathbb R \times \mathbb R^{k}} \left\{ \gamma - \mathbb P_n \left[ \Lambda(\gamma + \lambda^t \Phi_m(\theta,.))\right]\right\},
\end{equation}
where $\Lambda:\mathbb R \to \overline{\mathbb R}$ is a strictly convex, twice differentiable function such that $\Lambda'(0)=\Lambda''(0) = 1$ and $\Lambda(0)=0$. As previously, the existence of $\hat \theta_m$ requires the feasibility condition that the supremum of $\gamma - \mathbb P_n \left[ \Lambda(\gamma + \lambda^t \Phi_m(\theta,.))\right]$ is reached for a finite value of $(\gamma,\lambda) \in \mathbb R \times \mathbb R^{k}$, for at least one value of $\theta \in \Theta$. This condition relies essentially on the domain of $\Lambda$ being sufficiently widespread. We make the following additional assumptions. 
 
\begin{description}
\item A.6. The functions $x \mapsto \Vert \Phi(.,x) \Vert_\Theta$, $x \mapsto \Vert \nabla \Phi(.,x) \Vert_\mathcal N$ and $x \mapsto \Vert \Psi(.,x) \Vert_\mathcal N$ are dominated by a function $\kappa$ such that $\int \kappa^{4}(x) d \mu_0(x) < \infty$.
\item A.7. For all $x \in \mathcal X$ and for sufficiently large $m$, the map $\theta \mapsto \Phi_m(\theta,.)$ is twice continuously differentiable in $\mathcal N$ and we note $\nabla \Phi_m(\theta,.) = \partial \Phi_m(\theta,.)/\partial \theta$ and $\Psi_m(\theta,.) = \partial^2 \Phi_m(\theta,.)/\partial \theta \partial \theta^t$. 
\item A.8. The functions $ x \mapsto  \Vert \Phi_m(.,x) - \Phi(.,x) \Vert_\Theta$, $x \mapsto  \Vert \nabla \Phi_m(.,x) - \nabla \Phi(.,x) \Vert_\mathcal N$ and $x \mapsto  \Vert \Psi_m(.,x) - \Psi(.,x) \Vert_\mathcal N$ are dominated by a function $\kappa_m$ such that $\int \kappa_m^{4}(x) d \mu_0(x) = O(\varphi_m^{-4})$. 
\item A.9. The function $\Lambda''$ is bounded by a constant $K < \infty$.
\end{description}
Assumptions A.6 to A.8 are made to obtain a uniform control over $\Vert \hat \theta_m - \hat \theta \Vert$ for all $n \in \mathbb N$. The condition A.9 implies that $\Lambda$ is dominated by a quadratic function. In the MEM point of view, this condition is fulfilled for the log-Laplace transform $\Lambda_\nu$ of sub-Gaussian priors $\nu$.

\begin{theorem}[Robustness of GEL]\label{robustgel} If Assumptions 1 to 9 hold,
$$ n \Vert \hat \theta_m - \hat \theta \Vert^2 = O_P(n \varphi_m^{-2}) +  o_P(1).$$ 
Moreover, $\hat \theta_m$ is $\sqrt n$-consistent and asymptotically equivalent to the GEL estimator computed with exact constraint function $\Phi$ whenever $n \varphi_m^{-2}$ tends to zero.
\end{theorem}

By considering a situation with approximate operator, we extend the GEL model to a more general framework that gives a more realistic formulation of actual problems. The previous theorem gives an upper bound of the error caused by the use of the approximation $\Phi_m$ in place of the true function $\Phi$. By this result, we aim to provide an insight on convergence conditions that are necessary for asymptotic efficiency when dealing with an approximate operator. 

\section{Proofs} 
\subsection{Proof of Theorem \ref{gelmemth}}
\noindent Let $\mathcal S_\theta = \{ w \in \mathbb R^n: \mathbb P_n(w) \in \mathcal M_\theta \}$ and $\mathcal F_w = \{ \mu \in \mathcal P(\mathbb R^n): \mathbb E_\mu(W) = w\}$. We use that $ \inf_{\mu \in \mathcal F_w} \ \mathcal K(\mu, \nu_0) = \Lambda_{\nu_0}^*(w)$ (see \cite{MR1429928}). Let $ \Pi(\mathcal M_\theta) = \left\{ \mu \in \mathcal P(\mathbb R^n), \ \mathbb E_\mu \left[ \mathbb P_n(W)\right] \in \mathcal M_\theta \right\}$, we have the equality
$$ \hat \theta = \arg \min_{\theta \in \Theta} \ \inf_{\mu \in \Pi(\mathcal M_\theta)} \ \mathcal K(\mu,\nu_0) = \arg \min_{\theta \in \Theta} \ \inf_{w \in \mathcal S_\theta} \ \inf_{\mu \in \mathcal F_w} \ \mathcal K(\mu,\nu_0),  $$
which can be written
$$ \hat \theta = \arg \min_{\theta \in \Theta} \ \inf_{w \in \mathcal S_\theta} \Lambda_{\nu_0}^*(w) = \arg \min_{\theta \in \Theta} \ \inf_{w \in \mathcal S_\theta} \ \sup_{\tau \in \mathbb R^n}  \{ \tau^t w - \Lambda_{\nu_0}(\tau) \}.    $$
The feasibility assumption warrants that the extrema are reached. Hence, using Sion's minimax Theorem, we find
$$ \hat \theta = \arg \min_{\theta \in \Theta} \ \sup_{\tau \in \mathbb R^n}  \ \inf_{w \in \mathcal S_\theta} \{ \tau^t w - \Lambda_{\nu_0}(\tau) \},    $$
We know that $w=(w_1,...,w_n)^t \in \mathcal S_\theta$ if and only if $\sum^n_{i=1} w_i =n $ and $\sum^n_{i=1} w_i \Phi(\theta,X_i) = 0$. Thus, for a fixed value of $\tau$, the map $w \mapsto \tau^t w - \Lambda_{\nu_0}(\tau)$ can be arbitrarily close to $- \infty$ on $\mathcal S_\theta$ whenever $\tau$ is not orthogonal to $\mathds 1$ and $ \Phi(\theta,X)$. As a result, we may assume that $\tau = \gamma \mathds 1 + \Phi(\theta,X) \lambda$ for some $(\gamma,\lambda) \in \mathbb R \times \mathbb R^k$ without loss of generality. In this case, the map $w \mapsto \tau^t w - \Lambda_{\nu_0}(\tau)$ is constant over $\mathcal S_\theta$, equal to $n \gamma - \Lambda_{\nu_0}(\gamma \mathds 1 + \Phi(\theta,X) \lambda)$, which ends the proof. If $\nu_0 = \nu^{\otimes n}$, then $\Lambda_{\nu_0}(w) = \sum^n_{i=1} \Lambda_\nu(w_i)$ and we conclude easily.

\subsection{Proof of Theorem \ref{robustgel}} 
\noindent The proof of the results relies mainly on the uniform law of large numbers, using that the set $ \left\{ \Vert \Phi_m(\theta,.) \Vert, \Vert \nabla \Phi_m(\theta,.) \Vert, \Vert \Psi_m(\theta,.) \Vert, \ \theta \in \Theta, m \in \mathbb N \right\}$ is a Glivenko-Cantelli class of functions, consequently to A.6 and A.8. For all $\theta \in \Theta$, $v \in \mathbb R^k$, $x \in \mathcal X$, let
\begin{eqnarray*} h_n(\theta,v) & = & \left( \begin{array}{c} \mathbb P_n \left[\Phi(\theta,.) \Lambda'(v^t \Phi(\theta,.))\right] \\  \mathbb P_n \left[ v^t \nabla \Phi^t(\theta,.) \Lambda'(v^t \Phi(\theta,.))\right] \end{array} \right)\\
h_{m,n}(\theta,v) & = & \left( \begin{array}{c} \mathbb P_n \left[ \Phi_m(\theta,.) \Lambda'(v^t \Phi_m(\theta,.))\right] \\  \mathbb P_n \left[ v^t \nabla \Phi_m^t(\theta,.) \Lambda'(v^t \Phi_m(\theta,.))\right] \end{array} \right).
\end{eqnarray*}
The pair $(\hat \theta_m, \hat v_m)$ (resp. $(\hat \theta, \hat v)$) is defined as the unique zero over $\Theta \times \mathbb R^k$ of $h_{m,n}$ (resp. $h_n$). The condition A.9 implies that there exists a constant $K>0$ such that $\Lambda'(s) \leq K s + 1$ for all $s \in \mathbb R$. Hence, using successively the mean value theorem and Cauchy-Schwarz's inequality, we show that the contrast function $h_{m,n}$ converges uniformly on every compact set toward $h_n$ as $m \to \infty$, which warrants the convergence of $(\hat \theta_m, \hat v_m)$ toward $(\hat \theta, \hat v)$. For all $v \in \mathbb R^k$, the application $\theta \mapsto \nabla h_{m,n}(\theta,v)$ is continuous in a neighborhood on $\theta^*_m$ for sufficiently large values of $m$ by the condition A.7, as explicit calculation gives
$$ \nabla h_{m,n}(\theta,v) = \left( \begin{array}{cc} A_{m,n}(\theta,v) & D_{m,n}(\theta,v) \\ D_{m,n}^t(\theta,v) & V_{m,n}(\theta,v) \end{array} \right),  $$
where
\begin{eqnarray*} A_{m,n}(\theta,v) \!\!\!  &= &  \!\!\!  \mathbb P_n \left[ \Psi_m(\theta,.) v \Lambda'( v^t \Phi_m(\theta,.)) \right. \\
& & \ \ \ \ \  \left. + \nabla \Phi_m(\theta,.) v \ v^t \nabla \Phi^t_m(\theta,.) \Lambda''( v^t \Phi_m(\theta,.))\right] \\
D_{m,n}(\theta,v) \!\!\! & = &  \!\!\! \mathbb P_n \left[ \nabla \Phi_m(\theta,.) \Lambda'( v^t \Phi_m(\theta,.)) + \! \nabla \Phi_m(\theta,.) v \Phi^t_m(\theta,.) \Lambda''( v^t \Phi_m(\theta,.))\right] \\
V_{m,n}(\theta,v) \!\!\! &  = & \!\!\!  \mathbb P_n \left[ \Phi_m(\theta,.) \Phi_m^t (\theta,.) \Lambda''( v^t \Phi_m(\theta,.))\right].
\end{eqnarray*}
We define in the same way $A_n(\theta,v)$, $D_n(\theta,v)$ and $V_n(\theta,v)$ by replacing $\Phi_m$ by $\Phi$ in the expressions above. Using Cauchy-Schwarz's inequality, A.8 ensures the uniform convergence of $\nabla h_{m,n}$ toward $\nabla h_n$ on every compact set at the rate $\varphi_m$. Note $\rho_n$ the smallest eigenvalue of $\nabla h_n(\hat \theta, \hat v)$, we know from Theorem 3.2 in \cite{MR2031017} that $\mathbb P( \rho_n > \eta ) =$ \scriptsize{$O$}\normalsize $ (n^{-1})$ for sufficiently small $\eta >0$, since A.5 ensures that the limit of $\nabla h_n(\hat \theta, \hat v)$ as $n \to \infty$
is positive definite. Thus, for $c >0$ sufficiently small, consider the event $\Omega = \{ \rho_n > c \}$. Writing the Taylor expansion
$$ h_{m,n}(\hat \theta, \hat v) = \nabla h_{m,n}(\hat \theta_m, \hat v_m) \left( \begin{array}{c} \hat \theta - \hat \theta_m \\ \hat v - \hat v_m \end{array} \right) + o (\Vert \hat \theta_m - \hat \theta \Vert), $$
we deduce that on $\Omega$,
$$ \left( \begin{array}{c} \hat \theta_m - \hat \theta \\ \hat v_m - \hat v \end{array} \right) = - \left[ \nabla h_{n}(\hat \theta, \hat v) \right]^{-1} h_{m,n}(\hat \theta, \hat v) + O_P(\varphi_m^{-1}). $$
The Schur complement formula gives in particular
\begin{eqnarray*} \hat \theta_m - \hat \theta & = & - \left[ \hat D_n \hat V_n^{-1} \hat D_n^t \right]^{-1} \hat D_n \hat V_n^{-1} \ \mathbb P_n [ \Phi_m(\hat \theta,.) \Lambda'(\hat v^t \Phi_m(\hat \theta,.))] \\ 
& & \ \ \ \ + \ O_P(\varphi_m^{-1}) +  o_P(n^{-1}), 
\end{eqnarray*}
where $\hat D_n =D_n(\hat \theta, \hat v)$ and $\hat V_n=V_n(\hat \theta, \hat v)$ and where we used that $\hat v = O_P(n^{-1})$ (see for instance Theorem 3.2 in \cite{MR2031017}). Thus, on the event $\Omega$, 
$$ \Vert \hat \theta_m - \hat \theta \Vert \leq c \left| \left| \mathbb P_n [ \Phi_m(\hat \theta,.) \Lambda'(\hat v^t \Phi_m(\hat \theta,.))]   \right|\right| + O_P(\varphi_m^{-1}) +  o_P(n^{-1}).  $$
By construction, $\mathbb P_n [ \Phi(\hat \theta,.) \Lambda'(\hat v^t \Phi(\hat \theta,.))] =0$, which yields
\begin{eqnarray*} & & \left| \left| \mathbb P_n [ \Phi_m(\hat \theta,.) \Lambda'(\hat v^t \Phi_m(\hat \theta,.))]  \right| \right| \\ 
& \leq &  \mathbb P_n \left[ \Vert(\Phi_m(\hat \theta,.) - \Phi(\hat \theta,.)) \Lambda'(\hat v^t \Phi_m(\hat \theta,.))\Vert \right. \\
& & \ \ \ \ \ \ \left. + \ \Vert \Phi(\hat \theta,.) [ \Lambda'(\hat v^t \Phi_m(\hat \theta,.) - \Lambda'(\hat v^t \Phi(\hat \theta,.)) ] \Vert \right]  \\
& \leq & K \Vert \hat v \Vert \ \mathbb P_n \left[ \Vert \Phi_m(\hat \theta,.) \Vert \ \Vert \Phi_m(\hat \theta,.) - \Phi(\hat \theta,.) \Vert \right] + \mathbb P_n \Vert \Phi_m(\hat \theta,.) - \Phi(\hat \theta,.) \Vert \\
& & \ + K \Vert \hat v \Vert \ \mathbb P_n \left[ \Vert \Phi(\hat \theta,.) \Vert \ \Vert \Phi_m(\hat \theta,.) - \Phi(\hat \theta,.) \Vert \right], 
\end{eqnarray*}
as a consequence of A.9. We conclude that $\Vert \hat \theta_m - \hat \theta \Vert^2 \mathds 1_\Omega = O_P(\varphi_m^{-2}) +  o_P(n^{-1})$ by the condition A.8. On the complement of $\Omega$, $\Vert \hat \theta_m - \hat \theta \Vert$ can be bounded by the diameter $\delta$ of $\Theta$, yielding $\Vert \hat \theta_m - \hat \theta \Vert \mathds 1_{\Omega^c} =  o_P(n^{-1})$, which ends the proof.

\bibliographystyle{alpha}
\bibliography{bibliotheseabrev}

\end{document}